\documentclass[a4paper]{amsart}
\usepackage{amssymb}
\usepackage{amsmath} 
\usepackage{amsfonts}
\usepackage{amsthm}
\usepackage{amscd}
\usepackage{latexsym}
\usepackage{epsfig}
\usepackage{color}
\usepackage{hyperref}
\usepackage{mathrsfs}
\usepackage{graphics}
\usepackage{graphicx}
\usepackage{MnSymbol}
\usepackage [english]{babel}

\newtheorem{thm}{Theorem}[section]

\newcommand{\dx}{\, dx}

\newcommand{\D}{\partial}

\newcommand{\R}{\mathbb{R}}

\newcommand{\dv}{\mathrm{div}\,}

\newcommand{\sqc}{\mathfrak{C}}

\definecolor{mypink}{RGB}{219, 48, 122}
\definecolor{myblue}{RGB}{0, 0, 122}
\definecolor{mycyan}{rgb}{0.0, 0.55, 0.55}

\date{}

\begin{document}



\title[On the 3D tropical climate model with damping] {A regularity
  criterion for a 3D tropical climate model with damping}

\author{Diego Berti${}^\ast$}
\author{Luca Bisconti${}^\ast$}
\author{Davide Catania${}^{\ast\ast}$}

\begin{abstract}
  In this paper we deal with the 3D tropical climate model
  with damping terms in the equation of the barotropic mode $u$ and in
  the equation of the first baroclinic mode $v$ of the velocity,
  and we establish a regularity criterion for this system thanks to
  which the local smooth solution $(u, v, \theta)$ can actually
  be extended globally in time.
\end{abstract}

\maketitle

{\small${}^\ast${\small Universit\`a degli Studi di Firenze, Dipartimento di
    Matematica e Informatica ``U.~Dini'', Viale Morgagni 67/a, 50134 Firenze (FI), Italia} }\smallskip\smallskip
  
 {\small ${}^{\ast\ast}$\small {Universit\`a eCampus, Facolt\`a di Ingegneria, Via Isimbardi\, 10,
    I-22060 Novedrate (CO), Italia}}

\section{Introduction} \label{sec:intro}
We consider the following 3D tropical climate model with damping, i.e.
\begin{equation}\label{TCM-gen}
  \begin{aligned}
    &\D_t u+ (u\cdot \nabla) u -\nu\Delta u + \sigma_1 |u|^{\alpha -1}u +\nabla \pi+\dv  (v \otimes v)=0,\\
    &\D_t v+ (u \cdot \nabla) v -\eta\Delta v + \sigma_2 |v|^{\beta -1}v + (v\cdot \nabla) u  +\nabla\theta =0,\\
    & \D_t \theta + (u \cdot \nabla) \theta  -\mu \Delta \theta  + \dv  v =0,\\
    &\dv  u = 0,\\
    & u (x, 0) =u_0,\,\, v (x, 0)=v_0,\,\, \theta (x, 0)=\theta_0,
  \end{aligned}
\end{equation}
where $x\in \R^3$, $t \geq 0$, $u = (u_1(x,t), u_2(x,t), u_3(x,t))$
and $v = (v_1(x, t), v_2(x, t), v_3(x,t))$ denote the barotropic mode
and the first baroclinic mode of the velocity, respectively,
$\pi=\pi(x,t)$ indicates the pressure, $\theta = \theta(x,t)$ the
temperature, and $\nu > 0$, $\eta > 0$ and $\mu>0$.  Here,
$\sigma_1, \sigma_2>0$ and $\alpha, \beta \geq 1$ are the damping
coefficients (further appropriate restrictions will be introduced later on).

When $\nu = \eta = \mu = 0$ {(and $\sigma_1=\sigma_2=0$)}, the above system gives the original
tropical climate model derived by Frierson, Majda and Pauluis
\cite{Frierson} (see also \cite{Maj}).  Instead, in the case of $\nu>0$, $\eta>0$ and
$\mu=0$, \eqref{TCM-gen} reduces to {the viscous version of the
same model that has been analyzed by Li and Titi~\cite{Li-Titi} {(see also \cite{Li-Titi-1}).}

Local existence of strong solutions to the considered 3D model, {without damping,} and with
$\nu=1$, $\eta=1$ and $\mu=0$, has been established by Ma, Jiang and
Wan in \cite{Ma-Jiang-Wan}.  Global well-posedness of solutions to a
tropical climate model with dissipation in the equation of the first
baroclinic mode of the velocity, under the hypotheses of small initial
data, was studied by Wan \cite{Wan}, Ma and Wan \cite{Ma-Wan}.
In fact, the issue of global regularity has been investigated in a
number of articles and (partially) addressed, by introducing
  suitable hypotheses on the initial data, modified viscosity and diffusivity, or
  by inserting damping terms in the equations for $u$, $\nu$, and $\theta$.
Some of these studies, whose contents are related to the present analysis,
  are briefly recalled here below.

For the 3D case with damping, Yuan and Chen \cite{Yuan-Chen-2021}
studied the global regularity of strong solutions assuming
$\sigma_1>0$ and $\alpha\geq 4$, and removing the damping term from
the equation $\eqref{TCM-gen}_2$, i.e. setting $\sigma_2=0$ (see also
Yuan and Zhang, in \cite{Yuan-Zhang}).  In
\cite{Yuan-Zhang}, the authors proved global regularity assuming that
one of the following three conditions holds true:
\begin{align*}
  \text{(i)}& \qquad \alpha, \beta \geq 4, \\
  \text{(ii)}& \qquad \frac{7}{2} \leq \alpha<4,\,\,\,\, \beta\geq
  \frac{5\alpha+7}{2\alpha},
  \,\,\,\, \varpi\geq \frac{7}{2\alpha-5},\\
  \text{(iii)}& \qquad 3<\alpha \leq \frac{7}{2},\,\,\,\, \beta, \varpi\geq
  \frac{7}{2\alpha-5},
\end{align*}
where  $\varpi>0$ refers to an extra damping term,
i.e. $\sigma_3|\theta|^{\varpi-1}\theta$, with $\sigma_3\geq 0$, inserted in the left-hand
side of $\eqref{TCM-gen}_3$.

In \cite{Chen-Yuan-Zhang} the authors analyze the $d$-dimensional
\eqref{TCM-gen}, $d=2,3$, without damping terms, with only the
standard dissipation (i.e. $-\eta\Delta v$) of the first baroclinic
model of the velocity (and substituting $-\nu\Delta u$ with $\nu u$ {as
well as $-\mu\Delta \theta$ with $\mu\theta$)} and choosing a special class of initial data
$(u_0,v_0,\theta_0)$ with $H^s$-norm arbitrarily large. Let us also
recall the recent papers \cite{Li-Zhai-Yin, Li-Deng-Shang} where the
authors prove well-posedness for the tropical climate model upon
selecting special classes of initial data.
 
We also mention an article of Zhu \cite{Zhu}, in which the 3D system
\eqref{TCM-gen} is considered with $\sigma_i=0$, $i=1,2,3$, and
fractional diffusion on the barotropic mode, with initial data in
$H^3(\R^3)$.  The author proves global existence of strong solutions
$(u, v, \theta )\in L^\infty(0, T; H^3(\R^3))$, for any $T>0$,
removing $-\eta\Delta v$ in $\eqref{TCM-gen}_2$, and replacing
$-\nu\Delta u$ with $\nu\Lambda^{2\chi} u = \nu (-\Delta^{1/2})^{2\chi} u$ in
$\eqref{TCM-gen}_1$, with $\chi\geq 5/2$.

In \cite{Wang-Zhang-Pan} the authors give a regularity criterion on the gradient of $u$ for the
system~\eqref{TCM-gen}, assuming $\sigma_i=0$, $i=1,2,3$, and initial data in $H^2(\R^3)$. 
A further regularity criterion, for the local-in-time smooth solution to the 3D
tropical climate model in the Morrey--Campanato space, is given in \cite{Wu-2019}. 

Analyses yielding criteria similar to the one developed in the present
work are carried out in \cite {Fan-Jiang} for the 3D MHD equations
(see also \cite{Alghamdi, GRZ, GR5}).  {We also mention some
regularity criteria for the 3D Boussineq equations having connection
with our analysis, i.e. \cite{G4, Wu-2021, Zhang1, Zhang2}.}

In the present paper we consider problem \eqref{TCM-gen} with
$3\le\alpha, \beta<4$, and $(u_0,v_0,\theta_0)\in H^s$, {$3/2<s\leq 2$,}
providing a regularity criterion to obtain the smoothness of the
solutions. {This case, which is rather far from verifying the
hypotheses used to prove global existence in the previous works,}
highlights somehow the relevance of the choice $\alpha \geq 4 $ in
order to obtain smooth solutions: without such a condition, the
introduction of suitable constraints seems to be necessary in order to
obtain regular global solutions.

In our main result (see Theorem~\ref{main}, below) we provide a
regularity criterion involving the barotropic mode and the first
baroclinic mode of the velocity in the homogeneous Besov space
${\dot B^0_{\infty, \infty}}$ (see, e.g. \cite{Kozono2000,
  Kozono2002}), that is: if
\begin{equation*}
  \int_{0}^{T} \big(\|u(t)\|^{\delta}_{\dot B^0_{\infty, \infty}} 
  + \|v(t)\|^{\gamma}_{\dot B^0_{\infty, \infty}}\big)\,dt<\infty,
\end{equation*}
with $\delta=\delta(\alpha)$ and $\gamma=\gamma(\beta)$ defined in
\eqref{gamma-delta} below, then the solution
{$(u, v, \theta)(t)$ can be smoothly extended after $T>0$.}

\section{Preliminaries and main result} \label{sec:preliminaries-main}
%
For $p \geq 1$, we indicate by $L^p=L^p(\R^n)$ the usual Lebesgue
space, endowed with norm $\|\, \cdot\, \|_p=\|\, \cdot\,
\|_{L^p}$, with moreover
$\|\, \cdot\, \|= \|\, \cdot\, \|_2 $, when $p=2$.  For $s>0$, we denote by
$W^{s,p}=W^{s,p}(\R^n)$ and
$\|\, \cdot\, \|_{s,p}= \|\, \cdot\, \|_{W^{s,p}}$ the Sobolev space and
its norm, respectively (see, e.g., \cite{Ad}). When $p=2$, we use the
notation $H^s=W^{s,2}$ and
$\|\, \cdot\, \|_{H^s}=\|\, \cdot\, \|_{s,2}$. Here $\dot H^s$ denotes
the standard homogeneous Sobolev space with norm
$\|\, \cdot\, \|_{\dot H^s}$. {In terms of the Fourier transform, homogeneous and inhomogeneous
Sobolev spaces can be written as follows
\begin{equation*}
  \dot{H}^s = \Big\{ f\, :\, \|f\|^2_{\dot{H}^s} =
  \int_{\R^3}|\xi|^{2s}|\hat f(\xi)|^2d\xi <\infty \Big\},
\end{equation*}
and
\begin{equation*}
  H^s = \Big\{ f\, :\, \|f\|^2_{{H}^s} =
  \int_{\R^3}(1+|\xi|^{2s})|\hat f(\xi)|^2d\xi <\infty\Big\}.
\end{equation*}
Note that $H^s$ is equal to $\dot{H}^s \cap L^2$. For $s>0$, we also
introduce the operator $\Lambda^s$, formally
defined as $\Lambda^s f=(-\Delta )^\frac{s}{2}f$, that is the Fourier
multiplier such that $\widehat{\Lambda^s f}(\xi)=|\xi|^{s}\hat{f}(\xi)$, for
$\xi \in \R^n$. Plainly, $\Lambda^2 f= -\Delta f$.

  Most of the estimates involving Besov and $BMO$ spaces, that we
  use in the following, have been established in \cite{Kozono2002}.
We refer to this paper for a detailed overview (see also 
\cite{Kozono2000, Kozono2003}) on the theory of Besov spaces
$ B^s_{p,q}= B^s_{p,q} (\R^n)$ (and homogeneous Besov spaces
$\dot B^s_{p,q}=\dot B^s_{p,q} (\R^n)$), with $0 < p, q \leq \infty$
and $s \in \R$, and also on the $BMO=BMO(\R^n)$ the space of functions of {\em bounded mean oscillation}. 
 
In the sequel we will use the symbols {$C$ to denote generic
constants,} which may change from line--to--line, but are not dependent
on the specific functions under consideration. 

\subsection{Some estimates} \label{subsec:utility-estimates}
Here below, we make explicit those key tools (from literature) which are instrumental in order to prove Theorem \ref{main}.

We need interpolation's inequalities for Sobolev spaces and the well-known Gagliardo--Nirenberg's inequality (see \cite{Nirenberg-1959}). 
Besides this estimate in the standard form, we also make use and recall two recent fractional versions of it:
as a consequence of \cite[Corollary 2.4]{Hajaiej-2011}, we have
{
\begin{equation}
\label{e:GN-Hajaiej}
\|\Lambda^r f\|_{q}\le C \|f\|_{p_1}^{1-\kappa}\|\Lambda^{s}f\|_{p_2}^\kappa,
\end{equation}
in the case $r, s \ge 0$, $1<q, p_1, p_2< \infty$ and $0\le \kappa \le 1$ satisfying the 
conditions
\begin{equation} \label{GN-conditions}
\frac{1}{q} = \left(
  \frac{1-\kappa}{p_1}+\frac{\kappa}{p_2}\right)-\frac{\kappa s-r}{n}, \,\, \mbox{ and } \,\, r \le \kappa s.
\end{equation}
In the following, when we refer to Gagliardo--Nirenberg's
  inequality, we mean formula \eqref{e:GN-Hajaiej}.

For the extremal case $q=\infty$, we use \cite[Corollary 1]{Mi-Bre-2019} in the form
\begin{equation}
\label{e:GN-Brezis}
\| f\|_{\infty}\le C \|f\|_{p_1}^{1-\kappa}\|\Lambda^{s}f\|_{p_2}^\kappa,
\end{equation}
where, as before,  $0= \left(
    \frac{1-\kappa}{p_1}+\frac{\kappa}{p_2}\right)-\frac{\kappa s}{n}.$

Notice that the above Gagliardo--Nirenberg-type estimates allow for further fractional generalizations, like those in \cite[Theorem~1]{Mi-Bre-2018} and \cite[Theorem~1]{Mi-Bre-2019}, along with their direct consequences.
%

In the sequel, we will also use the Kato--Ponce product estimate \cite{Kato-Ponce-1988} {(see also \cite{G-2019, Gu, KPV1991})}, i.e.
\begin{equation}\label{KP}
  \|\Lambda^s(fg)\|_p \leq C\big(\|f\|_{p_1} \|\Lambda^s g\|_{q_1}
  + \|g\|_{p_2}\|\Lambda^s f\|_{q_2} \big),
\end{equation}
with $s>0$, $1<p<\infty$,  $1< q_1, q_2< \infty$ and $1< p_1, p_2\leq
\infty$, such that
$\frac{1}{p} = \frac{1}{p_1} + \frac{1}{q_1} = \frac{1}{p_2} +
\frac{1}{q_2}$. We highlight that, in the right-hand side of \eqref{KP}, 
$L^\infty$-norms are allowed only for zero-order terms with respect to $\Lambda^s$.

Also, for $1<p<\infty$, we will use (see \cite{Kozono2000})
\begin{equation} \label{BMO} \|f \,\cdot\, g\|_p \leq
  C(\|f\|_p\|g\|_{BMO} + \|f\|_{BMO}\|g\|_p).
\end{equation}

  At last, we recall the following logarithmic inequality
(see, e.g., \cite[Lemma~1.4, pp. 4 and 6]{GR6}):
\begin{equation} \label{log-besov}
  \begin{aligned}
    \|f\|_{BMO} &\leq C\|f\|_{\dot B^0_{\infty, 2}}
    \\
    &\leq C \big(1+\|f\|_{\dot B^0_{\infty
        ,\infty}}\ln^{\frac{1}{2}}(1+\|f\|_{H^s})\big),
  \end{aligned}
\end{equation}
where, in particular, we used the embedding
$\dot B^0_{\infty, 2} \subset BMO$, in the case $s>n/2$ (see, e.g.,
\cite{Kozono2000, Kozono2002, Kozono2003}, for more details).

From now on, it will always be assumed $n=3$.

\subsection{Regularity result} \label{subsec:reg-result}
\begin{thm} \label{main} {Let
  $(u_0, v_0, \theta_0) \in H^s \times H^s \times H^s$, {$3/2< s \leq 2$}, 
  with $\dv  u_0 = 0$.} Assume that $3\le \alpha, \beta<4$.
  Let $(u,v, \theta)$ be a local solution of the system
  \eqref{TCM-gen}, defined on some time interval $[0,T)$, with
  $0 < T < \infty$, and having the following regularity
  \begin{equation} \label{reg-dichiata} u, v, \theta \in C([0, \tilde
    T]; H^s)\,\, \textrm{ and }\,\, u, v\in L^2(0, \tilde T; H^{s+1} ),
  \end{equation}
  for any $0<\tilde T <T$.  Then $(u,v,\theta)(t)$ can be extended
  beyond time $T$, with the same regularity as in
  \eqref{reg-dichiata}, and hence as a smooth solution, provided that
  \begin{equation} \label{data-Besov} 
    \int_{0}^{T} \big(\|u(t)\|^{\delta}_{\dot B^0_{\infty, \infty}} 
    + \|v(t)\|^{\gamma}_{\dot B^0_{\infty, \infty}}\big)\,dt<\infty,
  \end{equation}
  where
  \begin{equation}
  \label{gamma-delta}
  \delta=\frac{6(\alpha-1)}{3\alpha-5}>2 \,\, \mbox{ and } \,\, \gamma=\frac{6(\beta-1)}{3\beta-5} >2.
  \end{equation}
\end{thm}

\section{Proof of Theorem~\ref{main}}
\label{sec-stime}
The proof consists in proving suitable energy estimates for the
considered solution {$(u, v, \theta)(t)$} showing explicitly that it can be
extended after time $T>0$. Thus, the procedure is divided in a number
of steps in which we establish the needed bounds in $L^2$, $H^1$ and
$H^s$, with {$3/2<s\leq 2$.} These steps parallel the formal estimates in the
global existence results given in \cite{Yuan-Chen-2021} and
\cite{Yuan-Zhang}, although in our case they are carried out with
different techniques borrowed from \cite{Bis-2020, Bis-2021, Dong-Wu-2019, G4, Ye, Zhang1}.

\subsection{$L^2$-estimates} \label{subsec:L2}
Taking the $L^2$-inner product of $\eqref{TCM-gen}_1$,
$\eqref{TCM-gen}_2$ and $\eqref{TCM-gen}_3$ with $u$, $v$ and
$\theta$, respectively, adding them up and integrating with respect to
$t$, we reach
\begin{equation} \label{L2-control}
  \begin{aligned}
    \|(u, v, \theta)(t)\|^2 + 2 \int_0^t\Big(\nu \|\nabla u\|^2   +& \eta
    \| \nabla v \|^2 +  \mu \|\nabla \theta\|^2 \\
    & +  \sigma_1  \|u\|_{\alpha+1}^{\alpha+1} 
    + \sigma_2 \|v\|_{\beta+1}^{\beta+1}\Big) \,d\ell
    =\|(u_0, v_0, \theta_0)\|^2,
  \end{aligned}
\end{equation}
where we used the notation
\begin{equation*}
  \|(u, v, \theta)(t)\|^2\doteq \|u(t)\|^2 +\|v(t)\|^2+\|\theta(t)\|^2,\,\,
  { \textrm{with $0 \leq t <T$},}
\end{equation*}
which will be adapted to the case of higher order norms, and the
following identities have been exploited
\begin{gather*}
  \int_{\R^3} \dv (v \otimes v) \cdot u \dx +
  \int_{\R^3} (v \cdot \nabla) u \cdot v \dx = 0,\\
  \int_{\R^3} \nabla \theta\cdot v \dx + \int_{\R^3} \dv v
  \cdot\theta \dx = 0,\\
  \int_{\R^3}(u\cdot \nabla )u\cdot u \dx =0, \,\, \int_{\R^3}(u\cdot
  \nabla )v\cdot v \dx =0 \textrm{ and }\, \int_{\R^3}(u\cdot \nabla
  )\theta\cdot \theta \dx =0.
\end{gather*}

Thanks to \eqref{L2-control}, for any $0<t<T$, it follows that
$u, v, \theta \in L^\infty (0, t; L^2)\cap L^2(0, t; H^1)$.

\subsection{\texorpdfstring{$\dot{H}^1$}{H1}-estimates} \label{subsec:H1-estimates}
Multiplying $\eqref{TCM-gen}_1$ by $-\Delta u$, integrating by
parts, we obtain
\begin{equation} \label{u-eq}
  \begin{aligned}
    \frac{1}{2} \frac{d}{dt} \|\nabla u(t)\|^2 + & \nu\|\Delta
    u(t)\|^2
    -  \sigma_1\int_{\R^3} u|u|^{\alpha-1}\cdot \Delta u\, dx \\
    & = \int_{\R^3}(u\cdot \nabla u)\cdot \Delta u \dx
    +\int_{\R^3}\dv  (v\otimes v)\cdot \Delta u \dx .
  \end{aligned} \hspace{-0.5 cm}
\end{equation}

Multiplying $\eqref{TCM-gen}_2$ by $-\Delta v$, we obtain
\begin{equation} \label{v-eq}
  \begin{aligned}
    \frac{1}{2} \frac{d}{dt} \|\nabla v(t)\|^2 + &\eta\|\Delta v(t)\|^2
     -  \sigma_2\int_{\R^3} v|v|^{\beta-1}\cdot \Delta v\, dx \\
    & = \int_{\R^3} (u\cdot \nabla ) v\cdot \Delta v\dx + \int_{\R^3}
    (v \cdot \nabla) u \cdot \Delta v\dx
    \\
    &\quad + \int_{\R^3} \nabla \theta \cdot \Delta v \dx.
  \end{aligned}
\end{equation}

Taking the $L^2$-product of $\eqref{TCM-gen}_3$ with $-\Delta \theta$,
we find
\begin{equation*} \label{theta-eq}
  \frac{1}{2} \frac{d}{dt}
  \|\nabla\theta (t)\|^2 +\mu \|\Delta \theta\|^2 = \int_{\R^3}
  (u\cdot \nabla )\theta \cdot \Delta \theta \dx + \int_{\R^3} \dv  v
  \, \Delta \theta\, dx.
\end{equation*}

Using calculations similar to those in \cite{Yuan-Zhang, Yuan-Chen-2021}, adding
\eqref{u-eq} and \eqref{v-eq}, we have that
\begin{equation} \label{stima-u-v-H1}
  \begin{aligned}
    \frac{1}{2}\frac{d}{dt}\big(&\|\nabla (u,v,\theta)\|^2 \big)
    + \nu\|\Delta u\|^2 + \eta\|\Delta v\|^2 +\mu \|\Delta \theta\|^2+
    \sigma_1\||u|^{\frac{\alpha-1}{2}}\nabla u\|^2 
    \\
    &\qquad +
    \frac{4\sigma_1(\alpha-1)}{(\alpha+1)^2}\|\nabla
    |u|^{\frac{\alpha+1}{2}}\|^2
    + \sigma_2\||v|^{\frac{\beta-1}{2}}\nabla v\|^2 +
    \frac{4\sigma_2(\beta-1)}{(\beta+1)^2}\|\nabla
    |v|^{\frac{\beta+1}{2}}\|^2
    \\
    &\quad = \int_{\R^3}(u\cdot \nabla) u\cdot \Delta u \dx
    +\int_{\R^3}\dv  (v\otimes v)\cdot \Delta u \dx +
    \int_{\R^3} (u\cdot \nabla ) v\cdot \Delta v\dx
    \\
    &\qquad + \int_{\R^3} (v \cdot \nabla) u \cdot \Delta v\dx +
    \int_{\R^3} \nabla \theta \cdot \Delta v \dx 
    \\
    &\qquad  
    +\int_{\R^3}
  (u\cdot \nabla )\theta \cdot \Delta \theta \dx + \int_{\R^3} \dv  v
  \, \Delta \theta\, dx  \doteq
    \sum_{i=1}^7 J_i,
  \end{aligned}\hspace{-1 cm}
\end{equation}
{where
\begin{equation*}
  \|\nabla (u, v, \theta)(t)\|^2\doteq \|\nabla u(t)\|^2 +\|\nabla v(t)\|^2+\|\nabla \theta(t)\|^2,\,\,\,
  \textrm{ with }\,\,\, 0\leq t< T.
\end{equation*}
}

Let us use $w$ to represent $u$, $v$ or even $\theta$. By applying H\"older's, Young's and 	Gagliardo--Nirenberg's inequalities, we obtain
\begin{equation*}
	\begin{aligned}
		J & \doteq \int_{\R^3}|u||\nabla v| |\Delta w| \, dx
		\le C\|u\|_6\|\nabla v\|_3 \|\Delta w\|
		\le C\|u\|_6^2 \|\nabla v\|_3^2+\varepsilon \|\Delta w\|^2
		\\
		& \le  C\|u\|_6^2 \|\nabla v\|\|\Delta v\|+\varepsilon \|\Delta w\|^2
		 = C\||u|^2\|_3 \|\nabla v\|\|\Delta v\|+\varepsilon \|\Delta w\|^2
		\\
		&\le  C\||u|^2\|_3^2 \|\nabla v\|^2+\varepsilon \|\Delta v\|^2+\varepsilon \|\Delta w\|^2.
	\end{aligned}
\end{equation*}	

From \eqref{BMO}, the interpolation\rq{}s inequality
\begin{equation*} 
	\|u\|_3 \le C \|u\|^\frac{2(\alpha-2)}{3(\alpha-1)}\|u\|_{\alpha+1}^\frac{\alpha+1}{3(\alpha-1)}
\end{equation*}
and \eqref{log-besov}, we have then
\begin{equation*} 
	\begin{aligned}
		J &\le C\|u\|_3^2 \|u\|_{BMO}^2 \|\nabla v\|^2+\varepsilon \|\Delta v\|^2+\varepsilon \|\Delta w\|^2
		\\
		&\le C\|u\|^{\frac{4(\alpha-2)}{3(\alpha-1)}}\|u\|_{\alpha+1}^\frac{2(\alpha+1)}{3(\alpha-1)}
		\|u\|_{BMO}^2 \|\nabla v\|^2+\varepsilon \|\Delta v\|^2+\varepsilon \|\Delta w\|^2
		\\
		&\le C\|u\|_{\alpha+1}^\frac{2(\alpha+1)}{3(\alpha-1)} \big(1+\|u\|^2_{\dot B^0_{\infty, \infty}}
		\ln(e+\|u\|_{H^s})\big) \|\nabla v\|^2+\varepsilon \|\Delta v\|^2+\varepsilon \|\Delta w\|^2.
	\end{aligned}
\end{equation*}
{To estimate lower order terms we used 
	\eqref{L2-control}, that, in particular, provide an uniform
	bound in time on $\|(u,v,\theta)(t)\|$, i.e $\sup_{0\leq t< T
	}\|(u,v,\theta)(t)\| \leq\|(u_0, v_0, \theta_0)\|.$}
	
Finally, from direct manipulations, the fact that $\frac{2(\alpha+1)}{3(\alpha-1)}<\alpha+1$,
and Young\rq{}s inequality with exponents $\frac{3(\alpha-1)}{2}$ and  $\frac{3(\alpha-1)}{3\alpha -5}$, we get
\begin{equation*}
	\|u\|_{\alpha+1}^\frac{2(\alpha+1)}{3(\alpha-1)}  \|u\|^2_{\dot B^0_{\infty, \infty}}
	\le C \big(\|u\|_{\alpha+1}^{\alpha+1} + \|u\|_{\dot B^0_{\infty, \infty}}^{2\frac{3(\alpha-1)}{3\alpha -5}}\big),
\end{equation*}
and so
\begin{equation*} 
	\begin{aligned}
		J &\le  C\|\nabla v\|^2 (1+\|u\|_{\alpha+1}^{\alpha+1} +\|u\|_{\alpha+1}^\frac{2(\alpha+1)}{3(\alpha-1)}
		\|u\|^2_{\dot B^0_{\infty, \infty}}\ln(e+\|u\|_{H^s})) +\varepsilon \|\Delta v\|^2+\varepsilon \|\Delta w\|^2
		\\
		& \le C\|\nabla v\|^2 \Big(1+\|u\|_{\alpha+1}^{\alpha+1} +\big(\|u\|_{\alpha+1}^{\alpha+1}
		+ \|u\|^{\frac{6(\alpha-1)}{3\alpha-5}}_{\dot B^0_{\infty, \infty}}\big)\ln(e+\|u\|_{H^s})\Big) +\varepsilon \|\Delta v\|^2+\varepsilon \|\Delta w\|^2.
	\end{aligned} \hspace{-1 cm}
\end{equation*}

Proceeding as for $J$, we deduce:
  \begin{equation*}
	\begin{aligned}
		J_1 &\leq \int_{\R^3}|u||\nabla u| |\Delta u| \, dx \\
		&\leq C\|\nabla u\|^2 \Big(1+\|u\|_{\alpha+1}^{\alpha+1} +\big(\|u\|_{\alpha+1}^{\alpha+1}
		+ \|u\|^{\frac{6(\alpha-1)}{3\alpha-5}}_{\dot B^0_{\infty, \infty}}\big)\ln(e+\|u\|_{H^s})\Big) +2\varepsilon \|\Delta u\|^2,
	\end{aligned}
\end{equation*}
\begin{equation*} 
	\begin{aligned}
		J_2 &\le  \int_{\R^3}|v||\nabla v| |\Delta u| \, dx 
		\\
		& \le C\|\nabla v\|^2 \Big(1+\|v\|_{\beta+1}^{\beta+1} +\big(\|v\|_{\beta+1}^{\beta+1}
		+ \|v\|^{\frac{6(\beta-1)}{3\beta-5}}_{\dot B^0_{\infty, \infty}}\big)\ln(e+\|v\|_{H^s})\Big) +\varepsilon \|\Delta v\|^2+\varepsilon \|\Delta u\|^2,
	\end{aligned} \hspace{-1 cm}
\end{equation*}
  \begin{equation*}
	\begin{aligned}
		J_3  &\leq  \int_{\R^3}|u||\nabla v| |\Delta v|\, dx \\
		&\leq C\|\nabla v\|^2 \Big(1+\|u\|_{\alpha+1}^{\alpha+1} +\big(\|u\|_{\alpha+1}^{\alpha+1}
		+ \|u\|^{\frac{6(\alpha-1)}{3\alpha-5}}_{\dot B^0_{\infty, \infty}}\big)\ln(e+\|u\|_{H^s})\Big) +2\varepsilon \|\Delta v\|^2,
	\end{aligned}
\end{equation*}
\begin{equation*} 
	\begin{aligned}
		J_4 &\le  \int_{\R^3}|v||\nabla u| |\Delta v| \, dx 
		\\
		& \le C\|\nabla u\|^2 \Big(1+\|v\|_{\beta+1}^{\beta+1} +\big(\|v\|_{\beta+1}^{\beta+1}
		+ \|v\|^{\frac{6(\beta-1)}{3\beta-5}}_{\dot B^0_{\infty, \infty}}\big)\ln(e+\|v\|_{H^s})\Big) +\varepsilon \|\Delta u\|^2+\varepsilon \|\Delta v\|^2,
	\end{aligned} \hspace{-1 cm}
\end{equation*}
and
  \begin{equation*}
	\begin{aligned}
		J_6  &\leq  \int_{\R^3}|u||\nabla \theta| |\Delta \theta|\, dx \\
		&\leq C\|\nabla \theta\|^2 \Big(1+\|u\|_{\alpha+1}^{\alpha+1} +\big(\|u\|_{\alpha+1}^{\alpha+1}
		+ \|u\|^{\frac{6(\alpha-1)}{3\alpha-5}}_{\dot B^0_{\infty, \infty}}\big)\ln(e+\|u\|_{H^s})\Big) +2\varepsilon \|\Delta \theta\|^2.
	\end{aligned}
\end{equation*}

{Finally, we observe that
\begin{equation*}
  J_5 + J_7 =\int_{\R^3} \nabla \theta \cdot \Delta v \dx +\int_{\R^3}
  \dv  v\,  \Delta \theta \dx =0.
\end{equation*}
}

Plugging the above estimates into \eqref{stima-u-v-H1}, we obtain
\begin{equation*} \label{stima-u-v-H1-chiusa}
  \begin{aligned}
    \frac{1}{2}\frac{d}{dt}&\big(\|\nabla (u,v,\theta)\|^2\big)
    + {\frac{\nu}{2}}\|\Delta u\|^2 + {\frac{\eta}{2}}\|\Delta v\|^2 +{\frac{\mu}{2}}\|\Delta \theta\|^2 \\
    &+ 
   {\frac{\sigma_1}{2}}\||u|^{\frac{\alpha-1}{2}}\nabla u\|^2 +
    \frac{4\sigma_1(\alpha-1)}{(\alpha+1)^2}\|\nabla
    |u|^{\frac{\alpha+1}{2}}\|^2
    \\
    &   +{\sigma_2}\||v|^{\frac{\beta-1}{2}}\nabla v\|^2 +
    \frac{4\sigma_2(\beta-1)}{(\beta+1)^2}\|\nabla
    |v|^{\frac{\beta+1}{2}}\|^2\\
    &\quad \le C \|\nabla (u,v,\theta)\|^2 \Big(1+\mathcal{A}
    +\big(\mathcal{A}+ \|v\|^{\gamma}_{\dot B^0_{\infty, \infty}} 
    + \|u\|^{\delta}_{\dot B^0_{\infty, \infty}}\big)\times\\
    & \qquad \times\ln(e+\|u\|_{H^s}+\|v\|_{H^s})\Big) 
  \end{aligned}\hspace{-1 cm}
\end{equation*}
where
\begin{equation*}
  \mathcal{A}=\|u\|_{\alpha+1}^{\alpha+1}+\|v\|_{\beta+1}^{\beta+1},\quad
  \gamma= \frac{6(\beta-1)}{3\beta-5}, \quad \delta=\frac{6(\alpha-1)}{3\alpha-5}.
\end{equation*}

In particular, for any $0\le t_*\le t <T$, we set
\begin{equation} \label{eq:y}
y(t)\doteq\sup_{t_*\le \ell\le t}\big(\|u(\ell)\|_{H^s}+\|v(\ell)\|_{H^s}\big),
\end{equation}
and, by applying Gronwall\rq{}s inequality, for any $0\leq t_* \le t <T$, we deduce
\begin{equation} \label{H1-control-new}
\begin{aligned}
\|\nabla(u, v, \theta)& (t)\|^2 + \int_{t_*}^{t}\left(
    {\frac{\nu}{2}}\|\Delta u\|^2 + {\frac{\eta}{2}}\|\Delta v\|^2 + {\frac{\mu}{2}}\|\Delta \theta\|^2\right) d\ell   \\
    & + {\int_{t_*}^t\Big(  \frac{\sigma_1}{2}\||u|^{\frac{\alpha-1}{2}}\nabla u\|^2 +
    \frac{4\sigma_1(\alpha-1)}{(\alpha+1)^2}\|\nabla
    |u|^{\frac{\alpha+1}{2}}\|^2\Big)\,d\ell}\\
    &+{\int_{t_*}^t \Big({\sigma_2}\||v|^{\frac{\beta-1}{2}}\nabla v\|^2 +
    \frac{4\sigma_2(\beta-1)}{(\beta+1)^2}\|\nabla
    |v|^{\frac{\beta+1}{2}}\|^2\Big)\,d\ell}\\
& \le \|\nabla(u, v, \theta)(t_*)\|^2\,e^{C \int_{t_*}^{t}\big( 1+ \|u\|^\delta_{\dot B^0_{\infty, \infty}} 
+\|v\|^\gamma_{\dot B^0_{\infty, \infty}}\big) \ln(e+\|u\|_{H^s}+\|v\|_{H^s})d\ell}\\
&\le C_*\big(e+ y(t)\big)^{C\varepsilon},
\end{aligned}
\end{equation}
with {$C_\ast$} positive constant only depending on $\|\nabla(u, v,
\theta)(t_*)\|^2$,  { and $t_\ast$  is such that}
\begin{equation*}
\int_{t_*}^{t}
\big(\|u(s)\|^{\delta}_{\dot B^0_{\infty, \infty}} +
\|v(s)\|^{\gamma}_{\dot B^0_{\infty, \infty}}\big)\,d\ell\le
\varepsilon, 
\end{equation*}
and $\varepsilon>0$ can be taken arbitrary small, in correspondence of suitable values of $t_*$, { because of \eqref{data-Besov}.}

\subsection{\texorpdfstring{$\dot H^s$}{Hs}-estimates, with \texorpdfstring{$\boldsymbol{3/2<s\leq 2}$}{3/2<s <= 2}}
\label{subsec:Hs-estimates}
In this last subsection, thanks to the previous $\dot{H}^1$-estimates, we provide 
  a proper energy inequality (see relation
  \eqref{e:finalfinal} below), at the level of the $H^s$-norm, which
  finally allows us to conclude the proof of Theorem~\ref{main}.
To this end, we apply the operator $\Lambda^s$ to each of
$\eqref{TCM-gen}_1$, $\eqref{TCM-gen}_2$ and $\eqref{TCM-gen}_3$ and
multiply in $L^2$ the resulting equations, respectively, by
$\Lambda^s u$, $\Lambda^s v$ and $\Lambda^s \theta$.

Applying $\Lambda^s$ to $\eqref{TCM-gen}_1$, and
multiplying the resulting equation in $L^2$, by $\Lambda^s u$, gives
\begin{equation*} \label{Hs-stima-diff-iniziale}
  \begin{aligned}
    \frac{1}{2} \frac{d}{dt} \|\Lambda^s & u(t)\|^2 + \nu\|\Lambda^{s+
      1}u(t)\|^2 
    \\
    = & -\int_{\R^3}\Lambda^s\big( (u\cdot \nabla) u\big) \cdot
    \Lambda^su \dx - \int_{\R^3}\Lambda^s\big( (v\cdot \nabla) v\big)
    \cdot \Lambda^su \dx
    \\
    &\quad  - \int_{\R^3}\Lambda^s\big( \dv  v\, v\big) \cdot
    \Lambda^su \dx - \sigma_1\int_{\R^3}\Lambda^s(|u|^{\alpha-1}u)\cdot
    \Lambda^s u\doteq \sum_{i=1}^4 K_i .
  \end{aligned}
\end{equation*}

Before going further with estimating $K_i$, $i=1,2,3, 4$, in order
to keep the notation compact, we introduce a non-negative,
time-dependent quantity of utility ${\sqc}={\sqc}(t)$, 
 that can be expressed as
\begin{equation} \label{sqc-star}
	{\sqc(t) \doteq \hat C + \chi_1\|\nabla u(t)\|^{\varrho_1} + \chi_2\|\nabla v(t)\|^{\varrho_2}
	+ \chi_3\|\nabla \theta(t)\|^{\varrho_3},} \,\,\, \hat C, \chi_i, \varrho_i\geq 0,\,\,\, i=1,2,3. \hspace{-0.2 cm}
\end{equation}
Here the constant $\hat C\geq 0$ is used to control the lower-order
terms $\|u\|$, $\|v\|$ and $\|\theta\|$ (which are bounded as a
consequence of \eqref{L2-control}); $\mathfrak{C}$ can change at any occurrence, from line--to--line. Let us notice
that, if necessary, the exponents $\varrho_i$ can be made explicit in the various subsequent steps,
even if it is in fact irrelevant for the purposes of the proof.

Let us then start with $K_2$. We have, by means of \eqref{KP},
with $\frac{1}{p_i}+\frac{1}{q_i}=\frac{1}{2}$, $i=1,2$, the following
inequality
\begin{equation} \label{K2-new}
  \begin{aligned}
  K_2& \leq  \left|\int_{\R^3}\Lambda^{s-1}\big( (v\cdot \nabla)
   v\big)\cdot
    \Lambda^{s+1}u \dx  \right|\\
  & \leq \|\Lambda^{s-1}((v\cdot\nabla)v)\|\|\Lambda^{s+1}u\|
  \\
  &\le C\Big(\|v\|_{p_1}\|\Lambda^{s-1}\nabla v\|_{q_1}+
  \|\nabla v\|_{p_2}\|\Lambda^{s-1}v\|_{q_2} \Big)\|\Lambda^{s+1}u\|
  \\
  &\le C\Big(\|v\|_{p_1}^2\|\Lambda^{s-1}\nabla v\|_{q_1}^2+
  \|\nabla v\|_{p_2}^2\|\Lambda^{s-1}v\|_{q_2}^2\Big) +
  \varepsilon\|\Lambda^{s+1}u\|^2
  \\
&  = C\Big(\|v\|_6^2\|\Lambda^{s-1}\nabla v\|_3^2+
  \|\nabla v\|^2_{\infty}\|\Lambda^{s-1}v\|^2\Big) +
  \varepsilon\|\Lambda^{s+1}u\|^2
  \\
  & {\leq C\Big(\| v\|_{1,2}^2\| \Lambda^{s-1}\nabla v\|_{3}^2+ \|\nabla
    v\|^2_{\infty}\|\Lambda^{s-1}v\|^2\Big) +
    \varepsilon\|\Lambda^{s+1}u\|^2,}
    \end{aligned}
\end{equation}
with the choice $(p_1,q_1)=(6,3)$ and $(p_2,q_2)=(\infty, 2)$, and by using the embedding $H^1\hookrightarrow L^6$.

{ We now conclude the estimate for $K_2$.
Using Gagliardo--Nirenberg's inequality, we infer 
\begin{equation} \label{s-l3}
 \begin{aligned}
	\|\Lambda^{s-1}\nabla v\|_3 &\leq C \|\nabla 
	v\|^{1-\kappa}\|\Lambda^{s+1}v\|^\kappa  
\end{aligned}
\end{equation}
with
%
$\kappa = 1-\frac{1}{2s}\le \frac{3}{4}.$}

Moreover, using again Gagliardo--Nirenberg's and the Young's
inequalities, we have
\begin{equation}
  \begin{aligned} \label{s-1-linfty-bis}
    \|\Lambda^{s-1} v\| &\leq C \|v\|^{2-s}\|\nabla v\|^{s-1}\\
    &\leq  C (\|v\| +\|\nabla v\|) \leq \sqc.
  \end{aligned}
\end{equation}
For the term $\|\nabla v\|^2_{\infty}\|\Lambda^{s-1}v\|^2$,
by applying \eqref{e:GN-Brezis}, we
deduce  that
\begin{equation} \label{e:forK2}
\|\nabla v\|_{\infty} \le C \|\nabla v\|^{\frac{2s-3}{2s}}\|\Lambda^{s+1} v\|^{\frac{3}{2s}}. 
\end{equation}
As a consequence of \eqref{s-l3}, \eqref{s-1-linfty-bis} and
\eqref{e:forK2}, we have, for any $3/2< s\le 2$ the following relations
\begin{equation*}
  \|\Lambda^s v\|_3 \le {\sqc} \|\Lambda^{s+1}v\|^{\frac{2s-1}{2s}},
  \quad \textrm{ and } \quad
  \|\nabla v\|_{\infty} \le {\sqc} \|\Lambda^{s+1}v\|^{\frac{3}{2s}}.
\end{equation*}
      
Thanks to the above estimates, Young's inequality,
and the fact that $\frac{2s-1}{s}, \frac{3}{s}<2$, then we get
\begin{equation}\label{hatK2}
  \begin{aligned}
  K_2
  &\le {\sqc}\big(\|\Lambda^{s+1}v\|^{\frac{2s-1}{s}}+\|\Lambda^{s+1}v\|^{\frac{3}{s}}\big)+ \varepsilon\|\Lambda^{s+1}u\|^2  
  \\
  &\le  {\sqc}  + \varepsilon \|\Lambda^{s+1}(u,v)\|^2,
    \end{aligned}
  \end{equation}
  where
  \begin{equation*}
  \|\Lambda^{s+1}(u, v)(t)\|^2\doteq \|\Lambda^{s+1}u(t)\|^2 +\|\Lambda^{s+1}v(t)\|^2,\,\,\,
  \textrm{with }\,\,\, 0\leq t<T.
\end{equation*}
  
As far as $K_1$ and $K_3$ are concerned, observe that we can follow what it has been
already done for $K_2$. Indeed, for $K_1$, we have
\begin{equation*}
\begin{aligned}
K_1 &\le \left|\int_{\R^3} \Lambda^{s-1}((u\cdot \nabla)u)\cdot \Lambda^{s+1}u \dx \right|
\\
 &\le  {\sqc}  + \varepsilon \|\Lambda^{s+1}(u,v)\|^2,
\end{aligned}
\end{equation*}
where we used the same exact calculations as in \eqref{hatK2} once substituting $v$ with $u$.
As far as $K_3$ is concerned, we have
{\begin{equation*}
\begin{aligned}
K_3 & \le \left| \int_{\R^3} \Lambda^{s-1}(\dv  v \, v)\cdot \Lambda^{s+1}u \dx \right|
\\
&\le \| \Lambda^{s-1}(\dv  v \, v)\|\| \Lambda^{s+1}u\|
\\
&\le C\Big(\|v\|_{p_1}\|\Lambda^{s-1} \dv v\|_{q_1}+ \|\dv v\|_{p_2}\|\Lambda^{s-1}v\|_{q_2} \Big)\|\Lambda^{s+1}u\|
  \\
  &\le C\Big(\|v\|_{p_1}^2\|\Lambda^{s}v\|_{q_1}^2+ \|\nabla v\|_{p_2}^2
  \|\Lambda^{s-1}v\|_{q_2}^2\Big) + \varepsilon\|\Lambda^{s+1}u\|^2.
\end{aligned}
\end{equation*}
}
Hence, with the same arguments used above, we conclude
\begin{equation} \label{hatK1-3}
K_1 + K_3 \le {\sqc}  + \varepsilon \|\Lambda^{s+1}(u,v)\|^2.
\end{equation}

Lastly, consider $K_4$. We have
\begin{equation*}
\begin{aligned}
K_4 &\le \left|\int_{\R^3}\Lambda^s(|u|^{\alpha-1}u)\cdot \Lambda^s u\, dx\right| 
\\
&= \left|\int_{\R^3}\Lambda^{s-1}(|u|^{\alpha-1}u)\cdot \Lambda^{s+1} u\, dx\right|
\\
&\le C\|\Lambda^{s-1}(|u|^{\alpha-1}u)\|^2 + \varepsilon\|\Lambda^{s+1}u\|^2
\\
&\le  C\left(\|\Lambda^{s-1} |u|^{\alpha-1}\|^2_3\|u\|^2_6 + \|\Lambda^{s-1} u \|^2_3 
\||u|^{\alpha-1}\|^2_6\right)+ \varepsilon\|\Lambda^{s+1}u\|^2,
\end{aligned}
\end{equation*}
where we used \eqref{KP} on $\|\Lambda^{s-1}(|u|^{\alpha-1}u)\|$ with $f
=|u|^{\alpha-1}$, $g=u$, $p_1=p_2=3$ and $q_1=q_2=6$. 
From {Gagliardo--Nirenberg's inequality, we have
\begin{equation*}
  \begin{aligned}
    \|\Lambda^{s-1}|u|^{\alpha-1}\|_3^2 &\leq C\||u|^{\alpha-1}\|^{2(1-z)}\|\nabla |u|^{\alpha-1}\|_3^{2z}
    \\
    &= C\||u|^{\alpha-1}\|^{\frac{4(2-s)}{3}}\|\nabla |u|^{\alpha-1}\|_3^{\frac{2(2s-1)}{3}},
\end{aligned}
\end{equation*}
since $
\frac{1}{3} = \frac{s-1}{3} + \frac{1-z}{2}$ implies $z=\frac{2s-1}{3}$. The same estimate holds true with $u$ instead 
of $|u|^{\alpha-1}$. These relations, along with the embedding $H^1\hookrightarrow L^6$ and the fact that
$ \alpha-1 < 3$ (so that 
$\|u\|_{2(\alpha-1)}\leq C\|u\|^{1-\lambda}\|\nabla u\|^{\lambda}$ for a suitable $0<\lambda\le 1$), yield
\begin{equation*} 
\begin{aligned}
K_4 &\leq C\left(\||u|^{\alpha-1}\|^{\frac{4(2-s)}{3}}\|\nabla |u|^{\alpha-1}\|_3^{\frac{2(2s-1)}{3}} \|u\|_{1,2}^2 
+ \|\nabla u\|_3^{\frac{2(2s-1)}{3}} \||u|^{\alpha-1}\|_{1,2}^2 \right) + \varepsilon\|\Lambda^{s+1}u\|^2
\\
&\leq C\|u\|_{2(\alpha-1)}^{\frac{4(2-s)(\alpha-1)}{3}}\||u|^{\alpha-2}\nabla u\|_3^{\frac{2(2s-1)}{3}} \|u\|_{1,2}^2
+C \|\nabla u\|_3^{\frac{2(2s-1)}{3}} \|u\|_{2(\alpha-1)}^{2(\alpha-1)}\||u|^{\alpha-2}\nabla u\|^2 
+ \varepsilon\|\Lambda^{s+1}u\|^2
\\
& \leq C\|\nabla u\|^{\frac{4(2-s)(\alpha-1)\lambda}{3}}\||u|^{\alpha-2}\nabla u\|_3^{\frac{2(2s-1)}{3}} \|u\|_{1,2}^2
+ C \|\nabla u\|_3^{\frac{2(2s-1)}{3}} \|u\|_{1,2}^{2(\alpha-1)} \||u|^{\alpha-2}\nabla u\|^2 
+ \varepsilon\|\Lambda^{s+1}u\|^2,
\end{aligned}
\end{equation*}
 { and so} 
\begin{equation} \label{e:hatK4i}
\begin{aligned}
K_4 &\le \sqc \||u|^{\alpha-2}\nabla u\|_3^{\frac{2(2s-1)}{3}} + \sqc
\|\nabla u\|_3^{\frac{2(2s-1)}{3}}
\||u|^{\alpha-2}\nabla u\|^2 + \varepsilon\|\Lambda^{s+1}u\|^2 
\\
& \doteq K_{41}+K_{42} +\varepsilon\|\Lambda^{s+1}u\|^2.
\end{aligned}
\end{equation}

We now focus on the first two terms of the right-hand side of this inequality. 
  
{For $K_{41}$, based on Gagliardo--Nirenberg's and Young\rq{}s inequalities, we have}
\begin{equation}
\label{e:delta}
\begin{aligned}
\||u|^{\alpha-2}\nabla u\|_3^\frac{2(2s-1)}{3} &\le \||u|^{\alpha-2}\|_6^\frac{2(2s-1)}{3} \|\nabla u\|_6^\frac{2(2s-1)}{3}
\\
&\le C \||u|^{\alpha-2}\|_6^\frac{2(2s-1)}{3} \Big(\|\nabla u\|^\frac{2(s-1)(2s-1)}{3s}\|\Lambda^{s+1}u\|^\frac{2(2s-1)}{3s}\Big)
\\
&= C\|u\|_{6(\alpha-2)}^\frac{2(\alpha-2)(2s-1)}{3} \|\nabla u\|^\frac{2(s-1)(2s-1)}{3s}\|\Lambda^{s+1}u\|^\frac{2(2s-1)}{3s}
\\
&\le C\|u\|_6^\frac{(2s-1)(2(\alpha-2)s-\alpha+3)}{3s} \|\Lambda^{s+1}u\|^\frac{(\alpha-3)(2s-1)}{3s}\\
&\quad \times \|\nabla u\|^\frac{2(s-1)(2s-1)}{3s}\|\Lambda^{s+1}u\|^\frac{2(2s-1)}{3s}
\\
& {\leq C \| u\|_{1,2}^\frac{(2s-1)(2(\alpha-2)s-\alpha+3)}{3s}} \|\Lambda^{s+1}u\|^{\frac{\alpha-1}{3}\frac{2s-1}{s}}\\
&\quad \times \|\nabla u\|^\frac{2(s-1)(2s-1)}{3s}
\\
&  {= C\| u\|_{1,2}^{\frac{(2s-1)[(2(\alpha-2)s-\alpha+3)+2(s-1)]}{3s}}}\|\Lambda^{s+1}u\|^{\frac{\alpha-1}{3}\frac{2s-1}{s}}
\\
&= {\sqc} \|\Lambda^{s+1}u\|^{\frac{\alpha-1}{3}\frac{2s-1}{s}}.
\end{aligned}
\end{equation}
Indeed, in the second line in the above relation, we used 
\begin{equation*}
\frac{1}{6}=\left(\frac{1}{2}-\frac{s}{3}\right)\kappa +\frac{1}{2}-\frac{\kappa}{2} \,\, \iff\,\, \kappa =\frac{1}{s},
\end{equation*}
which implies that
\begin{equation} \label{hatK4-secondopezzo}
  \|\nabla u\|_6 \le C \|\nabla u\|^{1-\kappa}\|\Lambda^{s+1}u\|^\kappa=
  C \|\nabla u\|^{\frac{s-1}{s}}\|\Lambda^{s+1}u\|^{\frac{1}{s}}.
\end{equation}
The fourth line in \eqref{e:delta} follows from 
\begin{equation*}
  \frac1{6(\alpha-2)}=\left(\frac12-\frac{s+1}{3}\right)\kappa +\frac16-\frac{\kappa}{6}
  \,\, \iff  \,\, {\kappa = \frac{\alpha-3}{2(\alpha-2)s}},
\end{equation*}
and so
\begin{equation} \label{hatK4-primopezzo}
  \|u\|_{6(\alpha-2)}\le C\|u\|_6^{1-\kappa} \|\Lambda^{s+1}u\|^\kappa
  = C\|u\|_6^{\frac{2(\alpha-2)s-\alpha+3}{2(\alpha-2)s}} \|\Lambda^{s+1}u\|^\frac{\alpha-3}{2(\alpha-2)s}.
\end{equation}

Observe that, in the last line of \eqref{e:delta},
since
\begin{equation*}
  \frac{\alpha-1}{3}\frac{2s-1}{s} < \frac{2s-1}{s}=2-\frac1{s}< 2,
\end{equation*}
we conclude that the first term in the right-hand side of \eqref{e:hatK4i}, by using Young's inequality,
can be estimated as
\begin{equation}
\label{e:delta-2}
K_{41}\le \sqc + \varepsilon \|\Lambda^{s+1} u\|^2.
\end{equation}

{ Similarly, for $K_{42}$, by means of  \eqref{hatK4-secondopezzo} and \eqref{hatK4-primopezzo}, we get}
\begin{equation*}
\begin{aligned}
\|\nabla u\|_3^{\frac{2(2s-1)}{3}} \||u|^{\alpha-2}\nabla u\|^2 & 
\le C\|\nabla u\|_3^{\frac{2(2s-1)}{3}} \||u|^{\alpha-2}\|_6 ^2\|\nabla u\|_3^2
\\
&= C\|\nabla u\|_3^{\frac{4(s+1)}{3}} \||u|^{\alpha-2}\|_6 ^2
\\
& \le C\|\nabla u\|^{\frac{2(s+1)(2s-1)}{3s}} \|\Lambda^{s+1} u\|^{\frac{2(s+1)}{3s}} \||u|^{\alpha-2}\|_6^2
\\
&=C\|\nabla u\|^{\frac{2(s+1)(2s-1)}{3s}} \|\Lambda^{s+1} u\|^{\frac{2(s+1)}{3s}} \|u\|_{6(\alpha-2)}^{2(\alpha-2)}
\\
&{\le {\sqc}  \|\Lambda^{s+1}u\|^{\frac{2(s+1)}{3s}+\frac{\alpha-3}{s}}}
\\
& ={\sqc} \|\Lambda^{s+1}u\|^\frac{3\alpha-7+2s}{3s}.
\end{aligned}
\end{equation*}

Hence, from $\frac{3\alpha-7+2s}{3s} < \frac{5+2s}{3s}<\frac{16}{9} <2$ and Young's inequality, we deduce
\begin{equation}
 \label{e:hatK4iii}
 K_{42}\le {\sqc} + \varepsilon \|\Lambda^{s+1} u \|^2.
 \end{equation}

In conclusion, from \eqref{e:hatK4i}, by using \eqref{e:delta-2} and \eqref{e:hatK4iii}, and renaming $\varepsilon$, 
we obtain
\begin{equation}  \label{hatK4}
K_4  \le {\sqc} +\varepsilon\|\Lambda^{s+1}u\|^2.
\end{equation}

Putting together \eqref{hatK2}, \eqref{hatK1-3} and \eqref{hatK4} we deduce
\begin{equation}\label{Hs-u-finale}
    \frac{1}{2} \frac{d}{dt} \|\Lambda^s u(t)\|^2 + (\nu - \varepsilon)\|\Lambda^{s+
      1}u(t)\|^2 \le {\sqc}(t) + \varepsilon \|\Lambda^{s+1}v(t)\|^2.
\end{equation}
for $0\leq t<T$.

Now, applying the operator $\Lambda^s$ to $\eqref{TCM-gen}_2$, and
multiplying the resulting equation by $\Lambda^s v$  in $L^2$, we get
\begin{equation} \label{Hs-stima-diff-v}
  \begin{aligned}
    \frac{1}{2} \frac{d}{dt} \|\Lambda^s & v(t)\|^2 + \eta\|\Lambda^{s+
      1}v(t)\|^2\\
       & = - \sigma_2\int_{\R^3}\Lambda^s(|v|^{\beta-1}v)\cdot
    \Lambda^s v\, dx
    -\int_{\R^3}\Lambda^s\big( (u\cdot \nabla) v\big) \cdot
    \Lambda^sv\dx  
    \\
     & \quad - \int_{\R^3}\Lambda^s\big( (v\cdot \nabla) u\big)
    \cdot
    \Lambda^sv \dx  - \int_{\R^3}\Lambda^s (\nabla \theta) \cdot
    \Lambda^sv \dx.
      \end{aligned}
    \end{equation}
 Similarly, for $\eqref{TCM-gen}_3$, we have
\begin{equation} \label{Hs-stima-diff-theta}
  \begin{aligned}
    \frac{1}{2} \frac{d}{dt} \|\Lambda^s & \theta(t)\|^2 +\mu\|\Lambda^{s+
      1}\theta(t)\|^2 \\
  &=-\int_{\R^3}\Lambda^s((u\cdot \nabla)\theta) \Lambda^s \theta \,dx
      - \int_{\R^3}\Lambda^s(\mathrm{div}\,v)\Lambda^s \theta\,dx.
      \end{aligned}
\end{equation}
{Therefore, adding \eqref{Hs-stima-diff-v} and \eqref{Hs-stima-diff-theta}, and observing that
\begin{equation*}
  \int_{\R^3} \Lambda^s(\nabla \theta) \cdot \Lambda^s v \dx +\int_{\R^3}
 \Lambda^s (\dv v)\,  \Lambda^s \theta \dx =0,
\end{equation*}
we reach
\begin{equation} \label{Hs-stima-diff-v-theta}
  \begin{aligned}
  \frac{1}{2} \frac{d}{dt}  \|\Lambda^s & ( v, \theta) (t)\|^2 + \eta\|\Lambda^{s+
      1}v(t)\|^2 +\mu\|\Lambda^{s+
      1}\theta(t)\|^2 \\
     &  = - \sigma_2\int_{\R^3}\Lambda^s(|v|^{\beta-1}v)\cdot
    \Lambda^s v\, dx
    -\int_{\R^3}\Lambda^s\big( (u\cdot \nabla) v\big) \cdot
    \Lambda^sv\dx  
    \\
     &  \quad - \int_{\R^3}\Lambda^s\big( (v\cdot \nabla) u\big)
    \cdot
    \Lambda^sv \dx  -\int_{\R^3}\Lambda^s((u\cdot \nabla)\theta) \Lambda^s \theta \,dx.
  \end{aligned}
\end{equation}
}

The first term in the right-hand side of \eqref{Hs-stima-diff-v-theta} can be estimated as done for $K_4$. 
After replacing $u$ with $v$ and $\alpha$ with $\beta$, we infer
\begin{equation*}
\sigma_2\left|\int_{\R^3}\Lambda^s(|v|^{\beta-1}v)\cdot
    \Lambda^s v \, dx \right| \leq 
{\sqc}  + \varepsilon \|\Lambda^{s+1} v\|^2,
\end{equation*}
where ${\sqc}$ has been introduced in \eqref{sqc-star}.

Moreover, using relation \eqref{KP}, and arguing as in \eqref{K2-new}--to--\eqref{hatK2}, we have
\begin{equation*}
\begin{aligned}
\bigg|\int_{\R^3} \Lambda^s \big((u\cdot\, &\nabla)v\big)  \cdot \Lambda^s v\,dx \bigg|
\\
& \leq \|\Lambda^{s-1}((u\cdot\nabla)v)\|\|\Lambda^{s+1}v\|
  \\
  &\le C\Big(\|u\|_{p_1}\|\Lambda^{s-1}\nabla v\|_{q_1}+
  \|\nabla v\|_{p_2}\|\Lambda^{s-1}u\|_{q_2} \Big)\|\Lambda^{s+1}v\|
  \\
  &\le C\Big(\|v\|_{p_1}^2\|\Lambda^{s-1}\nabla v\|_{q_1}^2+
  \|\nabla v\|_{p_2}^2\|\Lambda^{s-1}u\|_{q_2}^2\Big) +
  \varepsilon\|\Lambda^{s+1}v\|^2
  \\
&  = C\Big(\|v\|_6^2\|\Lambda^{s-1}\nabla v\|_3^2+ {
  \|\nabla v\|^2 _{\infty}\|\Lambda^{s-1}u\|^2\Big) +
  \varepsilon\|\Lambda^{s+1}v\|^2}
  \\
  & {\le C\Big(\| v\|_{1,2}^2\| \Lambda^{s} v\|_{3}^2+\|\nabla
    v\|^2 _{\infty}\|\Lambda^{s-1}u\|^2\Big) +
    \varepsilon\|\Lambda^{s+1}v\|^2}
\\
&\leq {\sqc} + \varepsilon \|\Lambda^{s+1}(u,v)\|^2,
    \end{aligned}
\end{equation*}
and, with very similar calculations to the ones in the previous case, we also have
\begin{equation*}
\bigg|\int_{\R^3} \Lambda^s\big((v\cdot  \nabla)u\big)\cdot \Lambda^s v\,dx\bigg|
\leq {\sqc} + \varepsilon \|\Lambda^{s+1}(u,v)\|^2
\end{equation*}
and, for the last term on the right-hand side of \eqref{Hs-stima-diff-v-theta}, it holds that
\begin{equation*}
  \left|\int_{\R^3} \Lambda^s((u\cdot \nabla)\theta) \Lambda^s \theta \,dx \right|
  \leq  {\sqc} + \varepsilon \|\Lambda^{s+1}(u,\theta)\|^2.
\end{equation*}
Hence, inserting the above estimates in \eqref{Hs-stima-diff-v-theta}, we infer
\begin{equation} \label{Hs-stima-diff-v-theta-finale}
  \frac{1}{2} \frac{d}{dt} \|\Lambda^s (v, \theta) (t)\|^2 + (\eta-\varepsilon)\|\Lambda^{s+ 1}v(t)\|^2
  + (\mu-\varepsilon)\|\Lambda^{s+ 1}\theta(t)\|^2 
   \leq {\sqc}(t) + \varepsilon \|\Lambda^{s+1} u (t)\|,
\end{equation}
for $0\leq t<T$.

Hence, putting together \eqref{Hs-u-finale} and
\eqref{Hs-stima-diff-v-theta-finale}, for $\varepsilon >0$ small enough,
we finally obtain 
\begin{equation}
\label{e:finalfinal}
\frac1{2}\frac{d}{dt} \|\Lambda^s (u,v,\theta)\|^2 + (\nu-\varepsilon)\|\Lambda^{s+1}u\|^2+ (\eta-\varepsilon)\|\Lambda^{s+1}v\|^2
+ (\mu-\varepsilon)\|\Lambda^{s+1}\theta\|^2 
\\
\le  \mathfrak{C(t)},
\end{equation}
for $0\leq t<T$.

From the above relation, integrating in time $(t_*, t)$, we reach
\begin{equation*}
  \begin{aligned}
&e+\|\Lambda^s (u,v,\theta)(t)\|^2 + \int_{t_*}^t \left( (\nu-\varepsilon)\|\Lambda^{s+1}u\|^2+ (\eta-\varepsilon)\|\Lambda^{s+1}v\|^2
+ (\mu-\varepsilon)\|\Lambda^{s+1}\theta\|^2\right) \, d \ell  \\
&\qquad \leq e+\|\Lambda^s
  (u,v,\theta)(t_*)\|^2 + \int_{t_*}^{t}\mathfrak{C}(\ell)\, d\ell\, .
  \end{aligned} \hspace{-1.2 cm}
\end{equation*}
Recalling the estimate \eqref{H1-control-new}, the definition of $y=y(t)$ introduced in \eqref{eq:y},
and actually re-defining such a quantity as
\begin{equation*}
y(t)\doteq\sup_{t_*\le \ell\le t}\big(\|u(\ell)\|_{H^s}+\|v(\ell)\|_{H^s} +\|\theta(\ell)\|_{H^s}\big),
\end{equation*}
and using
\begin{equation*}
  \int_{t_*}^{t}\mathfrak{C}(\ell)\, d\ell\le C_\ast\int_{t_*}^t\big(e+
  y(t)\big)^{C\varepsilon}\, d \ell\leq C_\ast(t-t_*)\big(e+
  y(t)\big)^{C\varepsilon},
\end{equation*}
with $C_\ast$ positive constant only depending $\|\nabla(u, v, \theta)(t_*)\|^2$ and $t_\ast$, we obtain
\begin{equation} \label{dettaglio-tecnico}
e+\|\Lambda^s (u,v,\theta)(t)\|^2 
\\
\le e+ \|\Lambda^s (u,v,\theta)(t_*)\|^2 + C_\ast (t-t_*) (e+y(t))^{C\varepsilon},
\end{equation}
and taking the supremum over $0\le \ell \le t$ {(up to introducing an additional
positive constant on the right-hand side of the above inequality, we can exploit relation \eqref{H1-control-new}, on the
left-hand side of \eqref{dettaglio-tecnico}, to reconstruct the $H^s$-norm and also bound
the lower-order terms),} we have
\begin{equation*}
  \begin{aligned}
e+y(t)^2 
&\le  e+ \|\Lambda^s (u,v,\theta)(t_*)\|^2 + C_\ast(e+y(t)^2)^{C\varepsilon},
  \end{aligned} \hspace{-0.2 cm}
\end{equation*}
with $C_\ast$ large enough.
Now, taking $t_*$ sufficiently close to $t$ so that $\varepsilon$  is such
that $C \varepsilon <1$, we conclude that
\begin{equation*}
(e+y(t)^2 )^{1-C \varepsilon}\le e +\|\Lambda^s (u,v,\theta)(t_*)\|^2  + C_\ast < +\infty
\end{equation*}
and hence
\begin{equation*}
\sup_{t_* \le \ell \le t }\big(e+\|\Lambda^s (u,v,\theta)(\ell)\|^2\big) < +\infty, \ \mbox{ for every $t \in (0,T)$.}
\end{equation*}
Since $t$ is arbitrary in the interval $(0,T)$, 
this concludes the proof of Theorem~\ref{main}.

\medskip

\noindent\textbf{Acknowledgements.} The authors are members of the Gruppo
Nazionale per l'Analisi Mate\-ma\-tica, la Probabilit\`a e le loro
Applicazioni (GNAMPA) of the Istituto Nazionale di Alta Matematica
(INdAM).

\end{document}